\documentclass[12pt]{article}
\usepackage{amsmath, amssymb, amsfonts, amsthm}
\def\CC{{\rm \kern.24em \vrule width.02em height1.4ex
depth-.05ex \kern-.26em C}}

\def\TagOnRight

\def\AA{{\it I}\hskip-3pt{\tt A}}

\def\QQ{\rlap {\raise 0.4ex \hbox{$\scriptscriptstyle |$}}
  {\hskip -0.1em Q}}

\newcommand{\be}{\begin{equation}}
\newcommand{\ee}{\end{equation}}
\newcommand{\bea}{\begin{eqnarray}}
\newcommand{\eea}{\end{eqnarray}}
\newcommand{\Bea}{\begin{eqnarray*}}
\newcommand{\Eea}{\end{eqnarray*}}

\newcommand{\bi}{\begin{itemize}}
\newcommand{\ei}{\end{itemize}}

\newtheorem{Definition}{Definition}[section]
\newtheorem{Theorem}[Definition]{Theorem}
\newtheorem{Lemma}[Definition]{Lemma}

\newtheorem{Conjecture}[Definition]{Conjecture}

\newtheorem{Revised Question}[Definition]{Revised Question}

\theoremstyle{remark}

\begin{document}

\title{The fourth smallest Hamming weight in the code of the
projective plane over $\mathbb{Z}/p \mathbb{Z}$}
\author{Bhaskar Bagchi\\
Theoretical Statistics and Mathematics Unit\\
Indian Statistical Institute\\
Bangalore - 560 059, India.\\
email: bbagchi@isibang.ac.in}
\date{}

\maketitle

\abstract{Let $p$ be a prime and let $C_p$ denote the $p$-ary code of the projective plane over ${\mathbb Z}/p\mathbb{Z}$. It is well known that the minimum weight of non-zero words in $C_p$ is $p+1$, and  Chouinard  proved that, for $p \geq 3$, the second and third minimum weights are $2p$ and $2p+1$. In 2007, Fack et. al. determined, for $p\geq 5$,  all words of $C_p$ of these three weights. In this paper we recover all these results and also prove that, for $p \geq 5$, the fourth minimum weight of $C_p$ is $3p-3$. The problem of determining all words of weight $3p-3$ remains open.}

\section{Introduction}

Throughout this paper, $p$ is a prime number and $\pi_p$ denotes the projective plane over the field $\mathbb{F}_p := \mathbb{Z}/p\mathbb{Z}$. Let $P$ be the point set of this plane and consider the $\mathbb{F}_p$-vector space $V_p:= \mathbb{F}_p^P$ of all functions from $P$ into $\mathbb{F}_p$. For $w \in V_p$, its support is the set $\{x \in P : w (x) \neq 0\}$, and the Hamming weight $|w|$ of $w$ is the size (cardinality) of its support.

We identify any line $\ell$ of $\pi_p$ with its indicator function. Thus $\ell$ is not only a set of points, but also the function $\ell : P \rightarrow \mathbb{F}_p$ given by $\ell(x)=1$ if $x\in \ell$ and $\ell(x)=0$ otherwise. Thus the lines have a second life as vectors in $V_p$. Consider the linear subspace $C_p$ of $V_p$ generated by all the lines of $\pi_p$. $C_p$ is called the ($p$-ary) code of $\pi_p$, and its elements are called the words of this code. $C_p$ inherits the Hamming weight function from $V_p$.

The vector space $V_p$ is equipped with the standard non-degenerate bilinear form $\langle \cdot, \cdot \rangle$ given by $\langle w, w^\prime \rangle =\sum\limits_{x \in P} w (x) w^\prime(x) ~(w, w^\prime \in V_p)$. As usual, two vectors $w, w^\prime \in V_p$ are called orthogonal, and we write $w \perp w^\prime$, if $\langle w, w^\prime \rangle =0$. The orthocomplement of $C_p$ with respect to this form is denoted by $C^\perp_p$, and is called the dual code of $\pi_p$.

It is well known (see \cite{EF}, \cite{PJ}) that $C^\perp_p$ is a subcode of $C_p$, of codimension one. Hence $\text{dim} (C_p) = \binom{p+1}{2}+1, \text{dim}(C^\perp_p)=\binom{p+1}{2}$.

All this is a miniscule part of the vast literature on the application of coding theory to incidence systems. The monograph \cite{EF} by Assmus and Key is a comprehensive account of this connection. The monograph \cite{PJ} by Cameron and van Lint is much more selective, but worth reading.

Since the zero word occurs in any linear code, it is common to define the minimum weight of a code to be the minimum of the Hamming weights of its non-zero words. It is also well known (see \cite{EF}, \cite{PJ} again!) that the minimum weight of $C_p$ is $p+1$, and - up to multiplication by non-zero scalars - the words of this weight are the lines of $\pi_p$. In \cite{KL} Chuinard proved that the second minimum weight of $C_p$ is $2p$.   It of course follows that, for $p\geq 3$, the third minimum weight of $C_p$ is $2p+1$. In \cite{FF}, Fack et. al. proved that, for $p \geq 5$,  the only words of weight $\leq 2p+1$ in $C_p$ are the linear combinations of at most two lines. In this paper, we prove:-

\begin{Theorem} Let $p \geq 5$. Then, the fourth minimum weight of $C_p$ is $3p-3$. The only words of $C_p$ of Hamming weight $< 3p-3$ are the linear combinations of at most two lines of $\pi_p$.
\end{Theorem}

However, we are unable to determine all the words of $C_p$ (or even of $C_p^\perp$) of Hamming weight $=3p-3$. Since $C^\perp_p \subset C_p$, and since the words of weight $3p-3$ constructed in Theorem 5.2 of \cite{BB} (and reproduced in Lemma 2.2 below) are actually in $C^\perp_p$, Theorem 1.1 immediately implies:

\begin{Theorem}\label{1.2} Let $p\geq 5$. Then the second minimum weight of $C_p^\perp$ is $3p-3$ (and the minimum weight is, of course, $2p$).
\end{Theorem}

The result of Theorem \ref{1.2} was conjectured by this author in \cite{BB}. It was also conjectured there that the words given in Lemma 2.2 are the only words of weight $3p-3$ in $C_p^\perp$. This problem remains open.

Let's define the {\bf Hamming spectrum} of a code to be the set of all numbers which occur as weights of code words. Define a {\bf gap} in the Hamming spectrum to be an interval $[a,b] (a\leq b)$ such that the interval is disjoint from the spectrum, but $a-1$ and $b+1$ belong to the spectrum. In this language, Theorem 1.1 says that $[1,p], [p+2,2p-1]$ and $[2p+2,3p-4]$ are the first three gaps in the Hamming spectrum of $C_p$ for $p \geq 7$. This result, together with Lemma 2.7 below, leads us to make a very bold conjecture. Admittedly, we have almost no evidence in its favour.

\begin{Conjecture} ({\bf The spectral gap conjecture}) For $1 \leq k \leq \sqrt{p-2}$, the intervals $[kp+2,(k+1)(p+1-k)-1]$ are gaps in the Hamming spectrum of $C_p$.
\end{Conjecture}

The condition $k \leq \sqrt{p-2}$ here ensures that the right end points of these intervals are indeed at least as large as the corresponding left end points. There may well be other sporadic gaps within this bound, reflecting combinatorial complexities of $\pi_p$.

This paper is inspired by our previous construction of words of weight $3p-3$ and by the beautiful ideas in the paper [5] by Fack et. al. In particular, we systematise and explore their idea of using Moorhouse bases (see Lemma 2.6 below) to understand code words in $C_p$. This idea, used repeatedly in [5], leads to a key tool in this paper, namely Lemma 2.7 below. However, we make no use of their modified Moorhouse bases. This idea is replaced by Lemma 2.4, another key tool. We also point out that another main idea of [5] was the use of a lower bound due to Ball and Blokhuis on the size of a double blocking set in $\pi_p$. (A double blocking set in a projective plane is a set of points meeting every line in at least two points.) But we make no use of double blocking sets: they turn out to have been a distraction. Indeed, even the support $S$ of an arbitrary word $w$ of weight $=3p-3$ in $C_p$ cannot be a double blocking set. Using the arguments of Section 3, it can be seen that (when $p \geq 5$) either there are at least three lines disjoint from $S$ (if $w \in C_p^\perp$) or there are at least nine lines meeting $S$ in a single point (if $w \not\in C^\perp_p$).

The key tools on the geometry of arbitrary code words of $C_p$ are gathered together in Section 2. Finally, in Section 3 we combine these tools with elementary counting arguments to prove Theorem 1.1.
\section{General words of $C_p$}

The following lemma is well known.

\begin{Lemma} A word $w \in \mathbb{F}^{P}_p$ is in $C_p$ iff for all lines $\ell$ of $\pi_p, \sum\limits_{x\in \ell} w (x)$ is a constant independent of $\ell$.
\end{Lemma}

{\bf Proof:} The sum above is $\langle w, \ell \rangle$. Since the total number of lines as well as the number of lines through any point is $1(\text{mod} ~p)$, adding over all lines one sees that the alleged constant (if it exists) must be $\sum\limits_{x \in P} w (x)=\langle w, {\bf 1}\rangle$, where {\bf 1} is the constant function 1 on $P$. So the statement says that $w \in C_p \Leftrightarrow w \perp {\bf 1}-\ell$ for all lines $\ell$. Since these lines  span $C_p$, and $C_p^\perp$ is of codimension one in $C_p$, it follows that the set $\{{\bf 1}-\ell : \ell \text{ a ~line ~of } \pi_p\}$ spans $C^\perp_p$. Therefore the statement amounts to the trivial fact that $w \in C_p$ iff $w$ is orthogonal to $C_p^\perp$. \hfill${\Box}$

The next lemma is Theorem 5.2 of \cite{BB}. We reproduce it here for completeness. Note that, to prove the uniqueness assertion in this lemma, we use Lemma 3.4, proved much later. But there is no circularity since we shall not use the uniqueness part of Lemma 2.2 in what follows.

\begin{Lemma} Let $p$ be odd. Let $x$ be a point of $\pi_p$ and let $\ell_i$, $1 \leq i \leq 3$, be three distinct lines of $\pi_p$ through $x$. Also let $m$ be a line of $\pi_p$ not passing through $x$. Consider the point set $S=(\ell_1 \cup \ell_2 \cup \ell_3) \backslash (m \cup \{x\})$. Then there is a word in $C_p^\perp \subset C_p$ with support $S$ (and hence weight $=3p-3$). Up to multiplication by non-zero scalars, this is the unique word of $C_p$ with support $S$.
\end{Lemma}

{\bf Proof:} Set up homogeneous co-ordinates $[X:Y:Z]$ on $\pi_p$ such that the equation of $\ell_1,\ell_2,\ell_3$ is $X=0,Y=0$ and $X=Y$ (respectively), and the equation of $m$ is $Z=0$. It follows that $x$ is the point $[0:0:1]$ and $S$ consists of the points $[0:1:t],[1:0:t]$ and $[1:1:t]$, where $t$ varies over $\mathbb{F}^\ast_p := \mathbb{F}_p \backslash \{0\}$. Define a word $w : P \rightarrow \mathbb{F}_p$ by putting $w (\cdot)=0$ on $P \backslash S$, and defining $w$ on $S$ by $w ([0:1:t])=w ([1:0:t])=t$ and $w ([1:1:t])=-t$ for $t \in \mathbb{F}^\ast_p$.

We compute $\langle w, n \rangle =\sum\limits_{x \in n}w (x)$ for all lines $n$ of $\pi_p$. If $n=m$ or if $n$ is one of the $p-2$ lines through $x$ other than $\ell_1,\ell_2,\ell_3$, then $n$ is disjoint from $S$ and hence $\langle w,n\rangle =0$. If $n=\ell_i$ $(1 \leq i \leq 3)$ then $\langle w,n\rangle =\pm \sum\limits_{t \in \mathbb{F}^\ast_p}t=0$ since $p$ is odd.

Any other line $n$ has an equation of the form $Z=aX+bY$ where $(a,b)\neq (0,0)$. Then $n \cap S \subseteq \{[1:0:a],[0:1:b],[1:1:a+b]\}$. Therefore $\langle w,n\rangle =a+b-(a+b)=0$ in this case also. Thus $w$ is orthogonal to all lines. Hence $w \in C_p^\perp$. Clearly $|w |=\# (S)=3p-3$.

If $p=3$, then $S$ is the symmetric difference of two lines, and the uniqueness is fairly easy to prove. So assume $p \geq 5$. If $w^\prime$ is another word of $C_p$ with support $S$, then fix a point $u\in S$ and look at the word $w^{''}:= w^\prime - \frac{w^\prime(u)}{w(u)}w$. The support of $w^{''}$ is contained in $S\backslash \{u\}$, and so $|w^{''}| <3p-3$. Therefore, by Lemma 3.4 below, $w^{''}$ is an $\mathbb{F}_p$-linear combination of at most two lines. But the support of any non-zero word of this form meets some line in at least $p$ points, and $S$ does not share $p$ points with any line. So $w^{''}=0$. Thus $w^\prime =\frac{w^\prime (u)}{w(u)} w$, a non-zero scalar times $w$. This proves uniqueness. \hfill$\Box$

We now introduce

{\bf Notation}: For any two distinct points $x_1,x_2$ of $\pi_p$, let $x_1 \vee x_2$ denote the unique line of $\pi_p$ joining $x_1$ and $x_2$. Dually, for any two distinct lines $\ell_1,\ell_2$ of $\pi_p$, let $\ell_1 \wedge \ell_2$ denote the unique point of $\pi_p$ common to $\ell_1$ and $\ell_2$.

\begin{Definition} For any set $S$ of points of $\pi_p$, and a line $\ell$ of $\pi_p$, we shall say that $\ell$ is a {\bf passant}, {\bf tangent} or {\bf secant} to $S$ if $\# (S\cap \ell)=0,1$ or 2, respectively.
\end{Definition}

In terms of this notation and terminology, we have:

\begin{Lemma} Let $p$ be odd, and $S$ be the support of a word in $C_p$. For any point $x$ of $\pi_p$, define the subset $S_x$ of $S$ by
\begin{equation*}
S_x = \left\{ \begin{array}{l}
        \{y \in S: x \vee y \text{ is a tangent to } S\} \text{ if } x \not\in S, \\
        \{y \in S : y \neq x \text{ and } x \vee y \text{ is a secant to } S\} \text{ if } x \in S.
      \end{array} \right.
\end{equation*}
Then $S_x$ is contained in a line of $\pi_p$.

{\bf Proof:} It suffices to show that any three distinct points $y_1,y_2,y_3$ of $S_x$ are collinear. (The result is trivial if $\# (S_x) \leq 2$.) Suppose not. For $1 \leq i \leq 3$, let $\ell_i$ be the line $x \vee y_i$. Also, let $m$ be the line $y_2 \vee y_3$. Then $x,y_1 \not\in m$. Consider the set $S^\prime =(\ell_1 \cup \ell_2 \cup \ell_3) \backslash (m \cup \{x\})$. Let $w \in C_p$ be a word with support $S$ and let $w^\prime \in C_p^\perp$ be a word with support $S^\prime$ (which exists by Lemma 2.2). Thus $\langle w, w^\prime \rangle =0$. But we have $S\cap S^\prime =\{y_1\}$ and hence $\langle w, w^\prime \rangle =w (y_1) w^\prime (y_1)\neq 0$ since $y_1$ belongs to the supports of $w$ and $w^\prime$. Contradiction. \hfill$\Box$
\end{Lemma}

{\bf Remark}: (a) If $p=2$ and $S$ is the complement of a line in $\pi_p$ then $S$ supports a word of $C_p$ and it violates the conclusion of Lemma 2.4 for any $x\in S$. Thus the hypothesis $p>2$ in Lemma 2.4 is essential.

(b) If $q >2$ is a power of $p$ and $C_q$ is the $q$-ary code of the projective plane $\pi_q$ over $\mathbb{F}_q$, then the construction of Lemma 2.2 actually produces a word of weight $3q-3$ in $C_q^\perp$ with a support of the same description. Therefore Lemma 2.4 is correct for the supports of words of $C_q$ for prime powers $q>2$.

The next lemma is a special case of a general observation in coding theory linking codes to matroids. The minimal supports of any linear code are precisely the circuits of an associated matroid. Lemma 2.5 occurs (essentially) as Theorem 2 in [5].

\begin{Lemma} Let $S$ be a set of points of $\pi_p$. Let $C^\ast_p$ be the $p$-ary code of the dual projective plane $\pi^\ast_p$. Then $S$ contains the support of some non-zero word of $C_p$ iff $S$ intersects every basis $B^\ast$ of $C^\ast_p$ such that $B^\ast$ consists of points of $\pi_p$ (i.e., lines of $\pi_p^\ast$).
\end{Lemma}

{\bf Proof:} Let $X$ be the subsystem of $\pi_p$ induced on $P\backslash S$. Thus, the point set of $X$ is $P \backslash S$, the lines of $X$ are the intersections with $P \backslash S$ of the lines of $\pi_p$. So, the dual incidence system $X^\ast$ has all points of $\pi^\ast_p$ as points, and the lines of $X^\ast$ are the lines of $\pi^\ast_p$ in $P \backslash S$. Consider the $p$-ary code $C_p (X^\ast)$ of $X^\ast$. It is a subcode of $C^\ast_p$. Since $C_p(X^\ast)$ is, by definition, spanned by the lines of $\pi^\ast_p$ in $P \backslash S$, it follows that there is a basis $B^\ast$ of $C_p^\ast$, consisting of lines of $\pi^\ast_p$, such that $B^\ast$ is disjoint from $S$ if and only if $C_p(X^\ast)=C^\ast_p$, i.e., iff $\text{dim} (C_p(X^\ast)) =$ $\text{dim} (C_p^\ast)$.

Note that, for any finite incidence system $X$, \linebreak  $\text{dim} (C_p(X^\ast))=\text{dim}(C_p(X))$. (This is because the incidence matrices of $X$ and $X^\ast$ are transposes of each other. Hence they have the same rank over $\mathbb{F}_p$.) In particular, we have $\text{dim}(C_p^\ast)=\text{dim}(C_p)$. (This particular case also follows from the fact that $\pi^\ast_p$ is isomorphic to $\pi_p$.)

Thus we get that there is a basis $B^\ast$ of $C_p^\ast$ (as above) disjoint from $S$ iff $\text{dim}(C_p(X))=\text{dim}(C_p)$. Now consider the restriction map $\rho : C_p \rightarrow C_p(X)$ which sends any $w\in C_p$ to its restriction to $P\backslash S$. Clearly $\rho$ is onto. Its kernel consists of all words of $C_p$ with support contained in $S$. The rank-nullity theorem implies that $\text{dim}(C_p(X))=\text{dim}(C_p)$ iff the kernel of $\rho$ is trivial, i.e., iff $S$ does not contain the support of any non-zero word of $C_p$. \hfill${\Box}$

Let $\tilde \pi_p$ denote the affine plane over $\mathbb{F}_p$, and let $\tilde P$ be its point set. Let $\tilde C_p$ denote the $p$-ary code of $\tilde \pi_p$, i.e., the subcode of $\mathbb{F}_p^{\tilde P}$ generated by the lines of $\tilde \pi_p$. Moorhouse [6] gave a beautiful combinatorial construction of a basis of $\tilde C_p$. In [5], the authors observe without proof that this basis ``extends" to a basis of $C_p$ by adjoining a  single suitable line. In the following, we include the simple proof of this fact for the sake of completeness.

\begin{Lemma} Let $\ell$ be a line of $\pi_p$. Enumerate the points on $\ell$ as $x_i, ~0\leq i \leq p$. For $1 \leq i \leq p$, let $L_i$ be a set of $i$ lines $\neq \ell$ through $x_i$. Also, let $m$ be any line through $x_0$ ($m$ may or may not be equal to $\ell$). Put $B_0 = \bigcup\limits^p_{i=1} L_i$ and $B=\{m\} \cup B_0$. Then $B$ is a basis of $C_p$. (We shall refer to any such basis $B$ as a Moorhouse basis of $C_p$.)
\end{Lemma}

{\bf Proof:} Let $\tilde \pi_p$ be the affine plane obtained from $\pi_p$ by discarding the line $\ell$ and the points on $\ell$. Thus the point set of $\tilde \pi_p$ is $\tilde P =P\backslash \ell$. Consider the restriction map $\rho :C_p \rightarrow \tilde C_p$ which maps any $w \in C_p$ to its restriction to $\tilde P$. Clearly $\rho$ is onto. Also, the kernel of $\rho$ consists of all words of $C_p$ with support contained in $\ell$. Since $p+1= \#(\ell)$ is the minimum weight of $C_p$, it follows that the kernel of $\rho$ is one dimensional: it is generated by $\ell$. Thus $\text{dim}(C_p)-\text{dim}(\tilde C_p)=1$. By the Moorhouse construction, $\rho(B_0)$ is a basis of $\tilde C_p$. Therefore $B_0$ is a linearly independent set. Also, if $D$ is the subcode of $C_p$ spanned by $B_0$, then $\rho$ restricts to a vector space isomorphism from $D$ onto $\tilde C_p$. Therefore $\text{dim}(D)=\text{dim}(\tilde C_p)$. Thus $D$ is a subcode of $C_p$ of codimension one. Note that $x_0$ does not belong to any of the lines in $B_0$. Therefore $x_0$ does not belong to the support of any word of $D$. It follows that if $m$ is a line of $\pi_p$ through $x_0$ then $m \not\in D$. Since $D$ is of codimension one in $C_p$ and since $B_0$ is a basis of $D$, it follows that $B=\{m\} \cup B_0$ is a basis of $C_p$. \hfill${\Box}$

{\bf Remark:} (a) Since $B$ is a basis of $C_p$ and $\# (B) =\binom{p+1}{2}+1$, it follows that $\text{dim} (C_p)=\binom{p+1}{2}+1$. Since $\text{dim} (C_p) +\text{dim} (C_p^\perp)=p^2+p+1=\binom{p+1}{2} +\binom{p+1}{2}+1$, it follows that $\text{dim}(C^\perp_p)=\binom{p+1}{2}$. Since  $C_p \cap C_p^\perp$ is clearly a codimension one subspace of $C_p$ (it is spanned by the differences of pairs of lines), Lemma 2.6 yields an independent proof of the fact that $C_p^\perp$ is a subcode of codimension one in $C_p$.

(b) It is known that if $q=p^e, e >1$, then the dimension of the $q$-ary code $C_q$ of $\pi_q$ is $\binom{p+1}{2}^e+1 < \binom{q+1}{2}+1$. Therefore the Moorhouse construction has no chance of working when $q$ is a ``genuine" power. This is why Lemma 2.6 (and its consequences, such as Lemma 2.7 below) is very specific to the case of a prime field.

The next lemma systematises an argument used repeatedly in \cite{FF}.

\begin{Lemma} Let $S$ be the support of a non-zero word of $C_p$ and let $x$ be a point of $\pi_p$. Then\\
\begin{enumerate}
\item[(a)] If $x \not\in S$ then there is a number $k$ in the range $2 \leq k \leq p+1$ such that there are at least $k$ lines $\ell$ through $x$ satisfying $\# (\ell \cap S) \geq p+2-k$. It follows that $\# (S) \geq k (p+2-k)$.

\item[(b)] If $x\in S$ and $S$ does not contain any line through $x$, then there is a number $k$ in the range $3 \leq k \leq p+1$ such that there are at least $k$ lines $\ell$ through $x$ satisfying $\# (\ell \cap S) \geq p+3-k$. It follows that $\#(S) > k(p+2-k)$.
\end{enumerate}
\end{Lemma}

{\bf Proof:} Order the lines through $x$ as $\ell_i,~ 0 \leq i \leq p$, in such a way that $\#(\ell_0 \cap S) \geq \# (\ell_1 \cap S) \geq \cdots \geq \# (\ell_p \cap S)$. If $x \not\in S$, then - as $x\in \ell_0$ - we have $\ell_0 \nsubseteq S$. If $x \in S$, then, by our assumption, $\ell_0 \nsubseteq S$. Thus, in either case, $\# (\ell_0 \backslash S) \geq 1$. If we also had $\# (\ell_i \backslash (S \cup \{x\})) \geq i$ for all $i~(1 \leq i \leq p)$, then there would be a Moorhouse basis of $C^\ast_p$ disjoint from $S$, contradicting Lemma 2.5 (since $S$ supports a non-zero word of $C_p$). So there is a number $k, ~2 \leq k \leq p+1$, such that $\# (\ell_{k-1} \backslash (S\cup \{x\}) \leq k-2$. Since $\# (\ell_{k-1})=p+1$, this means that (i) when $x\not\in S, ~\# (\ell_{k-1} \cap S) \geq p+2-k$ and (ii) when $x\in S, ~\# (\ell_{k-1} \cap S) \geq p+3-k$. Because of our numbering of the $\ell$'s, it follows that, for $0 \leq i \leq k-1, ~\# (\ell_i \cap S) \geq p+2-k$ in the first case and $\# (\ell_i \cap S) \geq p+3-k$ in the second case. Note that, in the second case, $k=2$ is impossible since $S$ contains no $\ell_i$. In case (a), the $k$ lines found cover $\geq k (p+2-k)$ points of $S$, so $\# (S) \geq k(p+2-k)$. In case (b), these lines cover $\geq k (p+2-k)$ points of $S\backslash \{x\}$, so that $\# (S) > k(p+2-k)$ in this case. \hfill${\Box}$

\section{Words of small weight}

In this section, we establish a series of lemmas on words of $C_p$ of Hamming weight $\leq 3p-3$, leading to a proof of Theorem 1.1.

\begin{Lemma} Let $S$ be the support of a word $w$ in $C_p$, such that $|w| \leq 3p-3$. Suppose no line of $\pi_p$ meets $S$ in $\geq p$ points. Then each point of $S$ is in at most two secants to $S$.
\end{Lemma}

{\bf Proof:} This is vacuous if $p=2$ or if $S$ is empty. So, assume $p \geq 3$ and $S \neq \emptyset$. Suppose, if possible, some point $x\in S$ is in three secants $\ell_i, ~1 \leq i \leq 3$. Say $\ell_i =x \vee x_i$, where $x_i \in S\backslash \{x\}$ are three distinct points.

By Lemma 2.4, there is a line $m$ containing $x_1,x_2,x_3$. Since $\ell_i \cap S =\{x,x_i\}$, $\langle w,\ell_i\rangle =w(x)+w(x_i),~ 1 \leq i \leq 3$. Therefore, Lemma 2.1 implies $w(x_1)=w(x_2)=w(x_3)=\lambda$ (say). Consider the word $w^\prime :=w -\lambda m\in C_p$. Let $S^\prime$ be the support of $w^\prime$. Note that $x\in S^\prime \subseteq S \cup m$, and the three lines $\ell_i$ are tangents to $S^\prime$ through $x$. Since $x$ is in the support of $w^\prime, ~w^\prime \neq 0$. Since $S^\prime \subseteq S \cup m$ and no line meets $S$ in $\geq p$ points, it follows that no line is contained in $S^\prime$. ($m \nsubseteq S^\prime$ since $x_i \in m \backslash S^\prime$.) Therefore, Lemma 2.7 (b) applies to the pair $(x,S^\prime)$.

Let $k$ be a number corresponding to $(x,S^\prime)$ guaranteed by Lemma 2.7 (b).

Thus $3 \leq k \leq p+1$ and there are $k$ lines $\ell$ through $x$ such that $\# (\ell \cap S^\prime) \geq p+3-k \geq 2$. Since there are also three tangent lines to $S^\prime$ through $x$, it follows that $k+3 \leq p+1$. Thus $3 \leq k \leq p-2$. Since $S^\prime \subseteq S \cup m$ and $x \not\in m$, these $k$ lines meet $S$ in $\geq p+2-k \geq 4$ points each.
Together with the three secants through $x$ they cover $\geq 1+3+k(p+1-k)$ points of $S$. Thus $k(p+1-k) \leq \# (S)-4<3(p+1-3)$. So $k(p+1-k)<3(p+1-3)$. But this is impossible since $3 \leq k \leq p+1-3$. \hfill${\Box}$

\begin{Lemma} Let $S$ be the support of a word $w \in C_p$. Suppose $S$ is not a line of $\pi_p$. Then,
\begin{itemize}
\item[(a)] If $|w| < 3p-3$, then each point outside $S$ is in at most two tangents to $S$.
\item[(b)] If $|w| =3p-3$, then each point outside $S$ is in at most three tangents to $S$. Indeed, if $x \not\in S$ is in three tangents to $S$, then each of the remaining $p-2$ lines through $x$ meets $S$ in exactly three points.
\end{itemize}
\end{Lemma}

{\bf Proof:} This is trivial if $p=2$ or if $S=\emptyset$. So assume $p \geq 3$ and $S$ is non-empty. Suppose some point $x \not\in S$ is in three tangents $\ell_i, ~1\leq i \leq 3$, to $S$. Say $\ell_i = x \vee x_i$, where $x_i \in S$ are distinct points. By Lemma 2.4, there is a line $m$ containing $x_1,x_2,x_3$. Since $S \cap \ell_i =\{x_i\}$, we have $\langle w, \ell_i\rangle =w(x_i)$. Therefore, Lemma 2.1 implies that $w(x_1)=w(x_2)=w(x_3)=\lambda$ (say). Consider the word $w^\prime :=w-\lambda m$. Let $S^\prime$ be the support of $w^\prime$. Note that $x \not\in S^\prime$, $S^\prime \subseteq S\cup m$, and the three lines $\ell_i$ through $x$ are passants to $S^\prime$. Since, by assumption, $S$ is not a line, we have $w \neq \lambda m$, and so $w^\prime \neq 0$. Thus, Lemma 2.7 (a) is applicable to the pair $(x,S^\prime)$.

Let $k$ be a number corresponding to $(x,S^\prime)$ guaranteed by Lemma 2.7 (a). Apart from the $k$ lines through $x$ meeting $S^\prime$ in $\geq p+2-k >0$ points each, there are also three passants to $S^\prime$ through $x$. So $k+3 \leq p+1$. Thus $2 \leq k \leq p-2$. Since $S^\prime \subseteq S \cup m$, and each line through $x$ meets $m$ in one point, we get $k$ lines through $x$ meeting $S$ in $\geq p+1-k$ points. Since $\langle w, \ell_i \rangle =w (x_i) \neq 0$ (as $x_i \in S$), Lemma 2.1 shows that $\langle w, \ell \rangle \neq 0$ for all lines $\ell$. Hence all lines meet $S$. Thus, apart from the $k$ lines through $x$ meeting $S$ is $\geq p+1-k$ points, the remaining $p+1-k$ lines through $x$ meet $S$ in $\geq 1$ point each. Together, these lines cover $\geq k(p+1-k)+p+1-k=(k+1)(p+1-k)$ points of $S$. So $(k+1)(p+2-(k+1))\leq \# (S) \leq 3 (p+2-3)$. Since $3 \leq k+1 \leq p-1$, this forces that $|w|=\# (S) =3p-3$, and $k=2$ or $k=p-2$. (This proves part (a)).

We first rule out the possibility $k=2$. If $k=2$, then we have two lines through $x$ meeting $S$ in $\geq p-1$ points and the remaining $p-1$ lines through $x$ meet $S$ in $\geq 1$ point. Since $\# (S)=3p-3$, it follows that two lines through $x$ meet $S$ in exactly $p-1$ points each and the remaining $p-1$ lines through $x$ are tangents to $S$. Let these tangents be $\ell_i=x \vee x_i$ $(1\leq i \leq p-1)$ where $x_1,\ldots,x_{p-1}$ are $p-1$ distinct points in $S$. As before, we have $w(x_1)=\ldots=w(x_{p-1})=\lambda \neq 0$, and, by Lemma 2.4, $x_i \in m$ for $1 \leq i \leq p-1$. Hence the $p-1$ lines $\ell_i$ through $x$ are passants to the support $S^\prime$ of $w^\prime$. Since the $k=2$ other lines through $x$ meets $S^\prime$ in (at least, hence exactly) $p+2-k=p$ points, it follows that $S^\prime =(\ell_p \cup \ell_{p+1}) \backslash \{x\}$, where $\ell_p, \ell_{p+1}$ are these two lines. But it is easy to verify (using Lemma 2.1) that the only words of $C_p$ with this support are the non-zero scalar multiples of the word $\ell_p-\ell_{p+1}$. So $w^\prime =\mu(\ell_p -\ell_{p+1})$ and hence $w=\lambda m +\mu (\ell_p-\ell_{p+1})$ where $\lambda \neq 0, \mu \neq 0$. But then $|w|=3p-2$ or $3p-1$ (according as $\mu=\pm \lambda$ or not). Contradiction. So $k \neq 2$.

Thus $k=p-2$. So we have $p-2$ lines through $x$ meeting $S$ in $\geq 3$ points each, and the remaining three lines through $x$ are tangents to $S$. Since $\# (S)=3p-3$, it follows that the first $p-2$ lines meet $S$ in exactly three points each. This proves (b). \hfill${\Box}$

\begin{Lemma} Let $S$ be the support of a non-zero word of $C_p$ of weight $< 3p-3$. Suppose no line meets $S$ in $\geq p$ points. Let $e_0,e_1,e_2$ denote the total number of passants, tangents and secants to $S$, respectively. Then, (a) $e_0 \leq 1$, (b) $e_1\leq p+2$, (c) $e_1+e_2 \leq \# (S)$.
\end{Lemma}

{\bf Proof:} Note that, if $x \not\in S$ is a point, then by Lemma 2.7 there are $k$ lines $(2 \leq k \leq p+1)$ through $x$ meeting $S$ in $\geq p+2-k$ points. Since $k(p+2-k)\leq \#(S) < 3 (p+2-3)$, we must have $k=2$, $p$ or $p+1$. Since no line meets $S$ in $\geq p$ points, $k\neq 2$. Thus $k=p$ or $p+1$. So there is at most one passant to $S$ through $x \not\in S$. But, if $\ell_1 \neq \ell_2$ were two passants, $x=\ell_1 \wedge \ell_2$ would be in two passants. So $e_0 \leq 1$.

A similar argument shows that, for any point $x\in S$, the corresponding number $k$ is $p$ or $p+1$. So there is at most one tangent through any $x\in S$, and if there is a tangent through $x$ then ($k=p$, and hence) $x$ is on no secant.

 By Lemma 3.2, any point outside $S$ is on at most two tangents. So, no three of the tangents to $S$ are concurrent. So, if $\ell$ is a tangent, each of the $p+1$ points on $\ell$ is on at most one more tangent. So $e_1 \leq p+2$.

Let $S_1 =\{x \in S: x$ is on a tangent to $S$\}. Put $S_2 =S \backslash S_1$. Each tangent meets $S_1$ in one point and we have seen that each point of $S_1$ is on only one tangent. So $\#(S_1)=e_1$. We have also seen that each secant (is disjoint from $S_1$, and hence) meets $S_2$ in two points, and - by Lemma 3.1 - each point of $S_2$ is on at most two secants. Thus $\#(S)-e_1 =\# (S_2) \geq e_2$. Therefore $e_1+e_2 \leq \#(S)$. \hfill${\Box}$

\begin{Lemma} Let $w \in C_p$ be of weight $< 3p-3$. Then $w$ is an $\mathbb{F}_p$-linear combination of at most two lines. (Hence $|w|=0, p+1, 2p$ or $2p+1$.)
\end{Lemma}

{\bf Proof:} Let $S$ be the support of $w$. We may assume $w \neq 0$. First suppose $S$ meets no line in $\geq p$ points. For $i \geq 0$, let $e_i$ be the number of lines meeting $S$ in exactly $i$ points. Counting in two ways the total number of lines and the number of point-line pairs $(x,\ell)$ with $x \in \ell \cap S$, we get $\sum e_i =p^2+p+1, ~\sum i e_i =(p+1) \# (S)$. Therefore,
\begin{eqnarray*}
(p+1)\# (S) &\geq& e_1 +2e_2 +3 \sum\limits_{i \geq 3} e_i \\
&=& e_1+2 e_2 +3(p^2+p+1-e_0-e_1-e_2)\\
&=& 3 (p^2+p+1) -3e_0-e_1-(e_1+e_2) \\
&\geq& 3(p^2+p+1) -3-(p+2)-\#(S)\\
&=& 3p^2+2p-2-\#(S).
\end{eqnarray*}
Here, the second inequality is by Lemma 3.3. Thus $|w| =\# (S) \geq \frac{3p^2+2p-2}{p+2} > 3p-4$. This is a contradiction. So, for every non-zero word $w \in C_p$ with $|w| < 3p-3$, there is a line $\ell$ which meets the support $S$ of $w$ in $\geq p$ points.

Now we complete the proof by finite induction on $|w|$. The result is trivial if $|w|=0$. So let $|w| >0$. Let $S$ and $\ell$ be as above. Since $w$ maps $S\cap \ell$ into $\mathbb{F}^\ast_p$, the pigeonhole principle yields two points $x_1 \neq x_2$ in $S\cap \ell$ such that $w(x_1)=w(x_2)=\lambda$ (say). Then $w^\prime :=w-\lambda \ell \in C_p$ and $|w^\prime| < |w|$. So by induction hypothesis, $w^\prime$ is an $\mathbb{F}_p$-linear combination of at most two lines. So $w$ is an $\mathbb{F}_p$-linear combination of at most three lines. But it cannot be a linear combination of three lines with non-zero coefficients since all such words have weight $>3p-3$. So $w$ is a linear combination of at most two lines. \hfill${\Box}$

{\bf Proof of Theorem 1.1:} Since $p \geq 5$, we have $p+1 <2p < 2p+1 < 3p-3$. Note that the first three numbers are precisely the weights of non-zero words which are $\mathbb{F}_p$-linear combinations of two lines. By Lemma 3.4, these three are the only numbers $<3p-3$ which occur as weights of non-zero words of $C_p$. Also, by Lemma 2.2, the number $3p-3$ occurs as an weight in $C_p$. So these four numbers are the four smallest non-zero weights in $C_p$. \hfill{$\Box$}

\end{document}